\documentclass[graybox]{svmult}

\usepackage{mathptmx,amsmath, amssymb}   
\usepackage{helvet}         
\usepackage{courier}        
\usepackage{type1cm}        
\usepackage{makeidx}         
\usepackage{graphicx}        
\usepackage{multicol}        
\usepackage[bottom]{footmisc}

\makeindex             

\def\paral{/\kern-0.55ex/}
\def\parals_#1{/\kern-0.55ex/_{\!#1}}
\def\bparals_#1{\breve{/\kern-0.55ex/_{\!#1}}}
\def\n#1{|\kern-0.24em|\kern-0.24em|#1|\kern-0.24em|\kern-0.24em|}

\newcommand{\B}{{\bf \mathcal B}}

\newcommand{\EE}{{\mathbf E}}

\newcommand{\F}{{\mathcal F}}

\def\h{{\mathfrak h}}

\newcommand{\R}{{\mathbf R}}
\newcommand{\Rc}{{\mathcal R}}

\def\L{{\mathcal L}}
\def\1{{\mathbf 1}}

\def\V{\mathbb V}

\def\le{\leq}
\def\ge {\geq}

\def\s.t.{\mathop {\rm s.t.}}

\def\Ric{\mathop{\rm Ric}}

\def\Hess{\mathop{\rm Hess}}

\def\<{\langle}
\def\>{\rangle}

\def\Lip{\mathrm {Lip}}

\def\V{\mathbb V}

\def\Lip{{\mathop {\rm Lip}}}

\def\det{\mathop{\rm det}}
\def\le{\leq}
\def\ge{\geq}
\def\paral{/\kern-0.55ex/}
\def\parals_#1{/\kern-0.55ex/_{\!#1}}
\def\bparals_#1{\breve{/\kern-0.55ex/_{\!#1}}}
\def\n#1{|\kern-0.24em|\kern-0.24em|#1|\kern-0.24em|\kern-0.24em|}
\def\f{\frac}

\begin{document}

\title*{Doubly Damped Stochastic Parallel Translations and Hessian Formulas}

 \titlerunning{Doubly damped stochastic parallel translation} 
  \authorrunning{Xue-Mei Li}
\author{Xue-Mei Li}

\institute{Department of Mathematics, Imperial College London, London SW7 2AZ. \email{xue-mei.li@Imperial.ac.uk.}  We dedicate this paper to Michael R\"ockner on the occasion of his 60's birthday.}

%
%
\maketitle

\abstract*{We study the fundamental solutions of a second order parabolic equation with possibly a zero order.}

\abstract{We study the Hessian of the  solutions of time-independent Schr\"odinger equations, aiming to obtain as large a class as possible of complete Riemannian manifolds for which the estimate $C(\f 1 t +\f {d^2}{t^2})$ holds.  For this purpose we introduce the doubly damped stochastic parallel transport equation, study them and make exponential estimates on them,  deduce a second order Feynman-Kac formula  and obtain the desired estimates. Our aim here is  to explain the intuition,  the basic techniques,  and the formulas which might be useful in other studies.}

\noindent
{\it AMS subject classification.} 60Gxx, 60Hxx, 58J65, 58J70.\\
{\it Keywords.}  Heat kernels, weighted Laplacian, Schr\"odinger operators, Hessian formulas, Hessian estimates.

\section{Introduction}
The probability distribution of a Brownian motion or a Brownian bridge are reference  measures with which we make $L^2$ analysis on the space of continuous paths (the Wiener space) and its subspaces of the pinned paths.  On the Wiener space, these are Gaussian measures and are well understood. The theory of the probability distributions of Brownian motion and Brownian bridges on  more general manifolds  is less developed. These include elliptic and semi-elliptic diffusion in an Euclidean space, with non-constant coefficients. For example we would like to describe the tail behaviour of the measure,
but how do we describe the set of path far away? Instead, we measure the size of the tails by checking whether a Lipschitz continuous function $f$  is exponentially integrable, or whether $\EE( e^{cf^2})$ is finite for a constant $c$.
In fact, a theorem of Herbst states that  if a probability measure  $\mu$ on $\R^n$ satisfies  the following  logarithmic Sobolev inequality
$\EE (f^2 \log (f^2))\le c_0 \EE |\nabla f|^2$, then $\EE (e^{\epsilon f^2})\le e^{\f{\epsilon }{(1-c_0 \epsilon)}}$  for any smooth function  $f:\R^n\to \R$  with $\|f\|_{L^2 (\mu)}=1$ where  $\epsilon$ is any number in $(0, 1/c_0)$ and. Also if  a mean zero  function $f$  is Lipschitz continuous with $\|f\|_{\Lip}=1$, then 
$\EE (e^{\alpha f})\le e^{c_0\alpha^2}$  holds for any number $\alpha$,
On a path or loop space, similar results hold \cite{Aida-Masuda-Shigekawa}. There are many other applications of Logarithmic Sobolev inequalities, see e.g. M. Ledoux's Saint Flour notes \cite{Ledoux} and the reference therein.
 To obtain such functional inequalities we make use of estimates on the fundamental solutions of Kolmogorov equations.

Let  $M$ denote a connected  smooth manifold with a complete Riemannian metric $g$.  Denote by $(g^{ij})$  the inverse of the Riemmanian metric $g=(g_{ij})$.
There exists on $M$ a strong Markov process with Markov generator $\f 12 \Delta$ where $\Delta$ is the Laplace-Beltrami operator which
 in local coordinates takes the form
 $$\Delta f(x)=\f 1 {\sqrt {\det g(x)}} \partial_i \left(\sqrt {\det g}\, g^{ij} \partial_j f  \right)(x).$$
 This stochastic process is said to be a Brownian motion.  If $Z$ is a vector field we denote by $L_Z$  Lie differentiation in the direction of $Z$ so  for a real valued function~$f$, $L_{Z}f =df(Z)$. Observe that any second order elliptic differential operator is of the form $\f 12 \Delta+L_Z$ where $\Delta$ is the Laplacian for the Riemannian metric induced by the operator, so they are generators of Brownian motions (with possibly a non-zero drift).  In this article we are mainly concerned  with  gradient drifts,  $Z=2\nabla h$ where $h:M\to \R$ is a smooth function.  In coordinates,   the operators are of the form
$\sum_{i,j}a_{i,j}(x) \f {\partial^2} {\partial x_i \partial x_j} +\sum_k \f {\partial} {\partial x_k}$, which is (locally) elliptic and is in general not strictly elliptic.

Set $\Delta^h=\Delta+2L_{\nabla h}$, this is  called the Bismut-Witten or the weighted Laplacian.  With respect to the weighted volume measure $e^{2h}dx$,  $\Delta^h$ is like a  Laplacian.  In particular, $\Delta^h=-(d+\delta^h)^2$ where $\delta^h$ is the adjoint of $d$ on $L^2(e^{2h}dx)$. All three operators, $d, \delta^h$ and $ \Delta^h$, extend to acting on differential forms.  Then $d+\delta^h$ and all its powers are essentially self-adjoint on $C_K^\infty$, the space of smooth compactly supported differential forms, for the details see \cite{Li-thesis}.  The densities of the probability distributions of the weighted Brownian motion are  the weighted heat kernel. These are the fundamental solution to the  equation $\f {\partial } {\partial t} =\f 12 \Delta^h$. There is also a commutative relation with which one can obtain  gradient estimates for the weighted heat kernels  under conditions on the Ricci curvature (without involving their derivatives).

We introduce the notations. Let  $\Ric_x: T_xM\times T_xM\to \R$ denote the Ricci curvature and let $\Rc$ denote the curvature tensor.  Let $\Ric^{\sharp}_x: T_xM\to T_xM$ denote the linear map defined  by the relation:
$\<\Ric^{\sharp}_x(u),v\>=\Ric_x(u,v)$. One of the novelties is to introduce the
  symmetrised tensor $\Theta$, see  \cite{Hessian},     \begin{equation}\label{Theta}
\left\< \Theta(v_2) v_1, v_3\right\>=  \left( \nabla_{v_3}  {\Ric}^{\sharp}\right) (v_1, v_2) -  \left( \nabla_{v_1}  {\Ric}^{\sharp}\right) (v_3, v_2)
 -   \left( \nabla_{v_2}  {\Ric}^{\sharp}\right) (v_1, v_3),
 \end{equation}
 where $v_1, v_2, v_3 \in T_{x_0}M$, 
 and to impose growth conditions on a bilinear map $\Theta^h$  from $T_{x_0}M\times T_{x_0}M$ to $\R$ instead of imposing conditions on $|\nabla \Ric^\sharp|$.  The bilinear form   is defined by the formula
$$\Theta^h(v_2) (v_1)=\f 12 
\Theta(v_2, v_1)+  \nabla^2 (\nabla h)( v_2,v_1)
+ \Rc( \nabla h, v_2)(v_1), \qquad v_1, v_2\in T_{x_0}M.$$

We are particularly interested in adding a zero order potential and  consider  a time independent Schr\"odinger equation, which is a parabolic partial differential equation of the form
\begin{equation}\label{1}
\f {\partial u} {\partial t} =(
\f 12 \Delta +L_{\nabla h}+V) u,
\end{equation}
where $u: \R\times M\to \R$ is a real valued function. For simplicity the zero order potential function $V: M\to \R$ will be assumed to be bounded and H\"older continuous and so $\Delta +2L_{\nabla h}+V$ is essentially self-adjoint on $C_K^\infty\subset L^2(M, e^{2h}dx)$,  the space of smooth functions with compact supports.

Our objectives are to  obtain global estimates for the Schr\"odinger semi-group $e^{(\f 12 \Delta^h +V) t} f$ where $f$ is a bounded measurable function,  for its gradient, and for its second order derivatives in terms of the geometric data of the Riemannian manifolds.
We will be also interested in such estimates for its fundamental solutions, which we denote by 
$p^{h,V}(t, x,y,)$, or $P^V(t,x,y)$ if $h$ vanishes identically or $p^h(t,x,y)$ if $V$ vanishes identically, and $p(t,x,y)$ if both $h$ and $V$ vanish. Similar notation, with capital $P$, e.g. $P_t^{h,V}$,  will be used to denote the corresponding semi-groups.

The commutative relation we mentioned earlier is as follows: the differential $d$ and the semi-group $e^{\f 12 \Delta^h}$ commute on $C_K^\infty$,  and consequently $d e^{\f 12 \Delta^h}$ solves  the heat equation on differential 1-forms:
$\f {\partial } {\partial t}\phi=\f 12 \Delta^h \phi$. If $M=\R^n$ this equation on differential 1-forms is an equation on `vector-valued'  functions.
Let us denote by $ \nabla^{h,*}$  the adjoint of $\nabla$ on  $L^2(e^{2h}dx)$, then this equation becomes: 
\begin{equation}
\label{1.3}\f {\partial } {\partial t}\phi=\f 12 \nabla^{h,*}\nabla \phi  -\f 12 \phi( {\Ric}^{\sharp}-2\nabla \nabla h).
\end{equation}
To see this we observe that, if $\nabla^*$ is the adjoint of $\nabla$ on $L^2(dx)$, then there is the Weitzenbock formula $\Delta^h=-  \nabla^{*}\nabla \phi{-\Ric}^{\sharp}(\phi)+2L_{\nabla h} \phi$, $L_{\nabla h}\phi=\nabla \phi (\nabla h) d(\phi(\nabla h)+d\phi(\nabla h, \cdot)$ where $\iota_{\nabla h}$ denotes the interior product. Also, for a differential 1-form  $\phi$ we apply the identity
$$\nabla^{h,*} \phi=  \nabla^{*}\phi-2\iota_{\nabla h}\phi $$ 
to see that
$$ \nabla^{h,*}\nabla_{\cdot} \phi=  \nabla^{*}\nabla_{\cdot} \phi-2\iota _{\nabla h}\nabla_\cdot  \phi. $$   

Equation (\ref{1.3}) inspired the study of the damped stochastic parallel translation
$$W_t: T_{x_0}M\to T_{x_t(\omega)}M$$ along a path $x_t(\omega)$ which solves the stochastic damped parallel translation equation
\begin{equation}
\label{wt}
\f{DW_t}{dt}=-\f 12 {\Ric}_{x_t}^{\sharp}(W_t)+\nabla_{W_t}\nabla h, \quad W_0=Id
\end{equation}
Here $Id$ denotes the identity map on $T_{x_0}M$ and 
$$\f{DW_t}{dt}:=\parals_t(x_\cdot(\omega)) \f{d}{dt}\left(   \parals_t^{-1} (x_\cdot(\omega))W_t\right)$$
is the covariant derivative along $x_t(\omega)$ and $\parals_t(x_\cdot(\omega)) : T_{x_0}M\to T_{x_t(\omega)}M$ denotes the standard stochastic parallel translation and  $\parals_t^{-1}(x_\cdot(\omega)):T_{x_t(\omega)}M \to T_{x_0}M$ is its inverse.  

Stochastic parallel translations along the non-differentiable sample paths of a Brownian motion can be  constructed  by a stochastic differential equation on the orthonormal frame bundle. This goes back to J. Eells,  K. D. Elworthy and P. Malliavin,  earlier attempts go back to K. It\^o and M. Pinsky.   M. Emery and M. Arnaudon studied parallel translations along a general  semi-martingale \cite{Emery-book}. The damped parallel translation goes back to E. Airault \cite{Airault}. 
The damped stochastic parallel translation takes into accounts of the effect of the Ricci curvature along its path and unwind it, leading to the magic well known formula:
$d e^{\f 12 \Delta^h t} f(v)=\EE df(W_t(v))$
for $(x_t)$ a Brownian motion with the initial value $x_0$. This holds  for all compact manifolds and 
for  more general manifolds.

The global estimates we are after are of the form
\begin{equation}\label{Hessian-type}
\left| \nabla^2 p(t, x_0, y_0)\right|  \le C\left( \f 1t +\f {d^2(x_0, y_0)}{t^2}\right), \quad t\in (0,1], \; x,y\in M.
\end{equation}
Such estimates (for $h=0$, $V=0$ and for compact manifolds )  were obtained in \cite{Malliavin-Stroock} and were generalised to other types of manifolds we refer to the references in \cite{Hessian}. We should remark that adequate care must be taken when generalising estimates from compact manifolds to non-compact manifolds.  For example taking a localising sequence of stopping times may not come for free and  any technique involving differentiating  a stochastic flow with respect to its initial point will likely need the additional assumption the strong 1-completeness \cite{Li-flow}. 
See also \cite{CLY}, \cite{Li-Yau}, \cite{Sheu}, \cite{HsuEstimates}, and \cite{Norris-93}.

Our goal is to establish these estimates for as large a class of manifolds as possible
and extend them to the operators $\f 12 \Delta +L_{\nabla h}+V$. If both $V$ and $h$ vanish,
these estimates are relevant for the study of the space of continuous loops and pinned paths using the Brownian bridge measure, e.g. the probability measure induced by a  Brownian motion conditioned to return to a point $y_0$  at time $1$. Naturally the Brownian motions with the symmetric drift $\nabla h$, which we refer as an $h$-Brownian motion, are also  candidates for such studies. We recall, that the $h$-Brownian bridge is a Markov process $x_t$ on $[0,1)$ with the Markov generator 
$$\f 12 \Delta ^h+ \nabla \log p(1-t, x, y_0),$$
and $\lim_{t\to 1} x_t=y_0$ where $y_0$ is the terminal value.
Observe that the corresponding damped parallel translation would be
$$\f{DW_t}{dt}=-\f 12 {\Ric}_{x_t}^{\sharp}(W_t)+\nabla_{W_t}\nabla h+\nabla_{W_t}  \nabla \log p(1-t, x_t, y_0), \quad t<1.$$

Gradient estimates on the semi-group associated with  the Brownian bridge  will naturally involve the second order derivative of  
$\log p(1-t, x, y_0)$. It is clear that the small time asymptotics of the Hessian are relevant, and estimates of the type (\ref{Hessian-type}) appear to be essential for analysing the Brownian bridge measure and useful  for the $L^2$ analysis of loop spaces. In this paper we explain the main formulas and constructions  rom \cite{Hessian} that leads  to these estimates.


%
%
%
%
%
\section{Summary of Results}
The following is  summary of some results from  \cite{Hessian}.
\begin{enumerate}
\item [(a)]We extend estimate  (\ref{Hessian-type})  to a more general class of manifolds replacing the linear growth condition on $|\nabla \Ric|_{op}$ by a linear growth condition on $\Theta$,  $\Theta$ being a symmetrised tensor obtained from $\nabla \Ric^\sharp$ after taking into accounts of the effects of~$\nabla h$ defined (\ref{Theta}). 
 \item [(b)] Our proof is based on an elementary  Hessian formula which will then lead to an integration by parts type Hessian formula.  For these formulas we  introduce a doubly damped stochastic parallel transport  equation which is defined using $\Theta$. It is natural to call these solutions  `doubly damped stochastic parallel translations'. We denote the solutions by  $W_t^{(2)}$, see Lemma \ref{lemma-doubly-damped}. We have the formula:
$$ \Hess (P^h_tf)(v_2, v_1) =\EE \left[\nabla df(W_t(v_2), W_t(v_1))\right]+\EE \left[df \left(W^{(2)}_t(v_1, v_2)\right)\right].$$
 Here $W_t$ denotes the damped stochastic parallel translation defined by (\ref{wt}). The Second Order Feynman-Kac Formula which does not involve the derivative of $f$  is given in Theorem
\ref{second-order}.

\item [(c)]  Such estimates will be also extended to the symmetric operators $\f 12 \Delta^h$  and $\f 12 \Delta^h+V$. Both operators are essentially self-adjoint on $C_K^\infty$.
By unitary transformations the drift term and the potential term can be treated almost exchangeably, however the drift and the zero order term do behave differently. 
For example we assume that $h$ is smooth and pose no direct assumptions on its its growth at the infinity while  the zero order term $V$ is only H\"older continuous and is assumed to be bounded. We obtain a second order Feynman-Kac formula, see Theorem \ref{second-order},  the H\"older continuity of $V$ is needed and is used to offset  singularities in some of the integrals
of the formulas. The modified doubly damped equation involves $\Theta^h$ instead of $\Theta$ itself.

\item[(d)] These estimates are refined for a subclass of manifolds with a pole, for which we make use of and obtain some nice estimates in terms of the semi-classical Brownian bridges,  a more careful study of the semi-classical Brownian bridge measure can be found in \cite{Li-ibp-sc}. See also \cite{generalised-bridges} for generalised Brownian bridges and \cite{Li-Thompson} for gradient estimates.
\end{enumerate}


\section{Key Ingredient}
Let $X(e)$ be smooth vector field on $M$ given by an isometric embedding $\phi:M\to \R^m$ and so $X(e)$ is the gradient of the real valued function $\<\phi, e\>$ where $e\in \R^m$. If $\{e_i\}$ is an o.n.b. of $\R^m$, this induces a family of vector fields $X_i(x)$ where $X_i(x)=X(x)(e_i)$. 

Let $F_t(x, \omega)$ denote the solution to the gradient
 SDE $$dx_t=\sum_{i=1}^m X_i(x_t)\circ dB_t^i+\nabla h(x_t)dt=X(x_t)\circ dB_t+\nabla h(x_t)dt$$ where $\circ $ denotes Stratonovich integration and $B_t$, $B_t=(B_t^1, \dots, B_t^m)$, is an $\R^m$-valued Brownian motion on a filtered probability space with the usual assumptions. Then  $F_t(x_0)$, the solution with the initial value $x_0\in M$, is a Brownian motion with the initial value $x_0$.

If $x_t$ is a semi-martingale, the stochastic damped parallel translation $\parals_t(x_\cdot(\omega))$ along $x_t(\omega)$, which is also denoted by $\parals_t$,  allows us to bring a vector in the tangent space of a solution path at time $0$
to its tangent space at time $t$, to differentiate it there and to bring it back to time $t$ by the inverse 
parallel translation $\parals_t^{-1}$. If $(x_t)$  is a Brownian motion with the initial value $x_0$, the damped stochastic parallel translation $W_t$ along its sample paths, where
$$\f{DW_t}{dt}=-\f 12 {\Ric}_{x_t}^{\sharp}(W_t)+\nabla_{W_t}\nabla h, \quad W_0=Id,$$
 compensates the effect of the Ricci curvature in equation (\ref{1.3})  and unwind it, leading to the magic well known formula, 
$$d e^{\f 12 \Delta^h t} f(v)=\EE \left[df(W_t(v))\right],$$
which holds trivially for compact manifolds and for manifolds with $\Ric-2\Hess h$ bounded from below and
for more general manifolds.

 For the second order derivatives of the fundamental solution of the heat kernel, we ought to differentiate $W_t$ with respect to its initial data, i.e. we differentiate $\parals_t^{-1}(x_\cdot)W_t(x_0)$  which is a map from $M$ to the space of linear maps which we denote by
$\L(T_{x_0}M; T_{x_0}M)$.

 We introduce the doubly damped stochastic parallel translation equation, whose solution we call doubly damped stochastic parallel transports/translations. Unlike damped parallel translations, the doubly damped ones involve genuine stochastic integrals (unless  the curvature vanishes) and it is a challenge to obtain exponential estimates. We also recall that the damped parallel translations are conditional expectations of  the spatial derivative of the solution to the gradient SDE. The doubly damped ones are obtained by differentiating the damped parallel translations, followed by taking conditional expectations. The beauty of it is that it satisfies the doubly damped stochastic parallel translation equation:
\begin{equation}\label{Wt2}
\begin{aligned} {D} v_t
=&  (- \f 12 \Ric  +\Hess h)^\sharp  (v_t )\,dt\\
& +\f 12 
\Theta^h(W_t(v_2))( W_t(v_1)) dt
+\Rc( d\{x_t\}, W_t(v_2)  )W_t(v_1).
\end{aligned}
\end{equation}

 We also introduce the notation $d\{x_t\}$, by which we mean integration with respect to the martingale part of $\{x_t\}$, see \cite[sec. 4.1]{Elworthy-LeJan-Li-book2} for detail. This allows us to give statements on a h-Brownian motion that is independent of its representation as a solution to a specific stochastic differential equation. In particular we may use 
 any of the two canonical representations:  (1) $x_t=\pi(u_t)$ where $u_t$ is the solution to the canonical 
 SDE on the orthonormal frame bundle $$du_t=H(u_t)\circ dB_t+\h_{u_t}\left( \nabla h(\pi(u_t)\right)\;dt,$$
 where $\h_u$ denotes the horizontal lift map at a given  frame $u$,
 and $\pi: OM\to M$ takes a frame, a point of $OM$, to its base point. (2) $x_t$ is the solution of a gradient SDE.
 If $u_t$ is the solution to the canonical SDE on the orthonormal frame bundle, then
 $d\{x_t\}$ is interpreted as $u_t \,dB_t$. If $(x_t)$ is the solution to
  to a gradient SDE driven by $X$ then $d\{x_t\}$  is interpreted as  $X(x_t) \,dB_t$. 
 
  \begin{lemma}\label{lemma-doubly-damped}
Suppose that the gradient SDE  is strongly 1-complete and suppose that
$v_1, v_2\in T_{x_t}M$ and $x_0\in M$ and $(x_t)$ is the solution to the gradient SDE. Let $W_t^{(2)}(v_1, v_2)$ denote the solution to the following covariant differential equation (the doubly damped stochastic parallel translation equation):
\begin{equation*}
\begin{aligned} {D} v_t
=&  (- \f 12 \Ric  +\Hess h)^\sharp  (v_t )\,dt\\
& +\f 12 
\Theta^h(W_t(v_2))( W_t(v_1)) dt
+\Rc( d\{x_t\}, W_t(v_2)  )W_t(v_1),\\
v_0&=0.
\end{aligned}
\end{equation*}
Then $W_t^{(2)}(v_1, v_2)$ is the local conditional expectation of $\nabla_{v_2} W_t(v_1)$ with respect to the filtration $\F_t^{x_0}:=\sigma\{x_s: s\le t\}$.
If furthermore the latter is integrable, then
$$W_t^{(2)}(v_1, v_2)= \EE \left\{\nabla_{v_2} W_t(v_1)\; \Big |\; \F_t^{x_0}\right\}.$$
 \end{lemma}

The proof for the lemma consists of stochastic calculus involving  $\f {D}{ds}W_t(j(s))$ where  $j$ is a parallel field with  $j(0)=v_2$,  along the normalised geodesic $\gamma$ 
with the initial condition $x_0$ and the initial velocity $\dot \gamma(0)=v_1$. 
Observe also that
$$\nabla_{v_2} W_t(v_1)=\f {D}{ds}|_{s=0} W_t(j(s)).$$

 Strong 1-completeness of an SDE is a  concept that is weaker than strong completeness, by the latter we mean the existence of a global solution to the SDE which is continuous with respect to the initial value.  Let $p=1, 2, \dots, n$ where $n$ is the dimension of the manifold. Roughly speaking, an SDE is strongly $p$-complete if for a.s. every $\omega$, and for all $t$,  $F_t(x, \omega)$ is continuous with respect to the initial point $x$ when $x$ is restricted to a   sub-manifold of dimension $p$ (or to a smooth $C^1$ curve if $p=1$).  The first example of an SDE which is complete (i.e. its solution from any initial point has infinite life time) and which is not strongly complete was given by K.D. Elworthy, prior to which it was generally believed that the two problems are equivalent. The concept  of strong $p$-completeness was introduced in \cite{Li-thesis, Li-flow} where we also give examples of  strongly $p-1$-complete SDEs which are not strongly $p$-complete and $n>2$ and $p\le n$. 
 For $n=1$ completeness is equivalent to strong completeness, similarly  for $n=2$, strong completeness is equivalent to strong completeness.
In \cite{Li-Scheutzow},   a non-strongly complete SDE on $\R^2$ is given: it has
one single driving linear Brownian motion and is driven a smooth bounded driving vector field.  We emphasise that, due to the fact that the exit time of $F_t(x, \omega)$ from a geodesic ball (even one with smooth boundary) is not necessarily continuous with respect to the initial point $x$, and it is not trivial to solve the strong  1-completeness by localisation. The strong 1-completeness for gradient SDEs was specially studied in \cite{Hessian}. See also the books \cite{Elworthy-book} and \cite{Kunita}.

\begin{remark}
If the gradient SDE is strongly 1-complete, 
$s\mapsto |W_t (\dot\gamma(s))|$ is continuous in $L^1$ and $\EE |T_{\gamma(s)}F_t|$ is finite, we  know that for all $f\in BC^1$, $d(P_t^hf)(v_1)=\EE df(W_t(v_1))$, \cite{Hessian, Li-flow, Elworthy-Li}.
From this we see immediately that
$$|d(P_tf)|_{L_\infty} \le  |df|_\infty \; \EE\left( e^{\int_0^t\rho^h(x_s) ds}\right),$$
where $\rho^h(x)=\sup_{|v|=1, v\in T_xM}\{- \f 12\Ric(v,v)+ Hess(h)(v,v)\}$. A more relaxed condition for this to hold can be obtained, but most of the assumptions here will be needed later.
If $\rho^h$ is bounded by $-K$ then we see immediately on direction of the characterisation for the Ricci curvature  to be bounded below by $K$, by taking $h=0$ in the earlier estimate, 
$$|d(P_tf)|_{L_\infty} \le  |df|_\infty \;e^{-Kt}.$$

\end{remark}

 We give below the second order analogue. Denote by $T_{x_0}F_t(v_0)$ the derivative flow of $F_t(x)$, it solves the equation
 $$\begin{aligned}dV_t&=(\nabla X)_{x_t}(V_t) \circ dB_t+(\nabla^2  h)_{x_t}(V_t) dt,\\
 V_0&=v_0.
 \end{aligned}
$$
Useful moment estimates on the derivatives flows for non-compact manifolds can be found in \cite{Li-moment}. Recall that $j(s)$ is a parallel field along the geodesic $\gamma$ with $\dot \gamma(0)=v_1$.
with the initial value $j(0)=v_2$.
\begin{lemma}\label{second-diff-lemma-1}
Suppose that $\Ric-2 \Hess (h)$ is bounded from below and that the gradient SDE  is strongly 1-complete.
Suppose also the statements (a) and (b) below hold.
\begin{enumerate}
\item [(a)] for every $s$, $\EE |T_{\gamma(s)}F_t|$ and $\EE |\nabla _{TF_t(\gamma(s))}W_t|$ are finite.
\item [(b)] $s\mapsto \EE\{ \f D {ds}W_t(j(s))\big| \; \F_t^{\gamma(s)}\}$ is continuous in $L^1(\Omega)$;
\end{enumerate}
Then  for all $f\in BC^2$,
\begin{equation}\label{basic-formula}
\Hess (P^h_tf)(v_2, v_1) =\EE \left[\nabla df(W_t(v_2), W_t(v_1))\right]+\EE \left[df \left(W^{(2)}_t(v_1, v_2)\right)\right].
\end{equation}
\end{lemma}

 From this lemma we immediately obtain the following estimate:
 $$\left| \Hess (P^h_tf)\right| _{L_\infty}
\le | \nabla df|_\infty \;  \EE  \left( e^{2\int_0^t\rho^h(x_s) ds}\right)
+|df|_{L_\infty}\; \EE  \left| W^{(2)}_t\right|.
$$

It is clear that estimation on $\EE  \left| W^{(2)}_t\right| $ will be useful,
this is given in \cite{Hessian} and  which we do not include here. As we shall see,  to obtain
Hessian estimates of the form (\ref{Hessian-type}), we will need to obtain exponential integrability of $ \left| W^{(2)}_t\right| ^2$,  such estimates will be given shortly after we explain why this is so.

The following are the basic assumptions.
 [{\bf  \underline{C1}.}]
\begin{itemize}
\item [(a)]  $\Ric-2 \Hess(h)\ge -K$; \item[(b)]  $\sup_{s\le t}\EE( \|W_s^{(2)}\|^2)<\infty$;
\item[(c)] for all  $f\in BC^2(M;\R)$, $v_1, v_2 \in T_{x_0}M$,  the elementary Hessian formula (\ref{basic-formula}) holds.
\end{itemize}
 
For the Schr\"odinger equation we first set, for $r<t$ and $x_0\in M$,
\begin{equation}
\V_{t-r, t}=(V  (x_{t-r})-V(x_0))e^{-\int_{t-r}^t  [V(x_s)-V(x_0)]ds}.
\end{equation}
Set also, \begin{equation}
N_{t} =\f 4 {t^2} \int_{\f t 2}^t \<X(x_s) dB_s, W_s(v_1)\>\int_0^{\f t 2} \<X(x_s) dB_s, W_s(v_2)\>.
\end{equation}

\begin{theorem}[Second Order Feynman-Kac Formula]
\label{second-order}
Suppose that {\bf  \underline{C1}} holds. Let  $V$ be a bounded H\"{o}lder continuous function.  Then for any $f\in \B_b(M;\R)$,
\begin{equation}\label{hessian}\begin{aligned}
 \Hess P_t^{h,V}f (v_1, v_2)
=&  e^{-V(x_0)t}\EE\left[ f(x_t)N_{ t}\right]
+ e^{-V(x_0)t}\EE \left[ f(x_t)\f 2 t\int_0^{t/2} \<X(x_s) dB_s,  W_{s}^{(2)}(v_1, v_2)\> \right]\\
&\quad+ e^{-V(x_0)t}\int_0^t  
\EE \left[f (x_t)\f {2 \V_{t-r, t} } {t-r}\int_0^{(t-r)/2} \<X(x_s) dB_s, W_{s}^{(2)}(v_1,v_2)\> \right] dr\\
&\quad+e^{-V(x_0)t} \int_0^t  \EE\left[ f (x_t)  \V_{t-r, t}N_{ t-r}\right]dr.
\end{aligned}
\end{equation}
\end{theorem}

For $h=V=0$, a version of the Hessian formula was first given in \cite{Elworthy-Li}
followed by another  in \cite{ArnaudonPlankThalmaier}.
A  version of  the Hessian formula for $h\equiv 0$ and $V\not =0$ was also given in \cite{Elworthy-Li}, however no proof was given. The doubly damped stochastic parallel translation equations were not present in either papers, nor were any extensive estimates given.   Hessian formula and estimates for non-linear potential, on linear space, were given  in \cite{Li-Zhao}. A formula for the Laplacian of
the semigroup $P_tf$  can be found in \cite{Elworthy-Li-form}.

 \begin{corollary} We assume $V(x_0)=0$. Then,
$$\begin{aligned}&\Hess p^{h,V}(t, x_0,y)\\
&=\Hess p_t^h(x_0,y)+
\int_0^t \int_M V(z) \Hess p^h(t-r, x_0,z) p_r^{h} (z,y)\EE[ e^{-\int_0^r V(Y_s^{r,z,y})}] dz dr,\end{aligned}$$
where $Y_s^{r,z,y}$ is the $h$-Brownian bridge with terminal value $r$, initial value $z$ and terminal value $y$.
\end{corollary}

Finally we indicate how to obtain estimates from these formulas. Let us take $V=0$ for simplicity, so the formula reads:
\begin{equation}\begin{aligned}
 \Hess P_t^{h}f (v_1, v_2)
=& \EE\left[ f(x_t)N_{ t}\right]
+ \EE \left[ f(x_t)\f 2 t \int_0^{t/2} \<X_i(x_s) dB_s,  W_{s}^{(2)}(v_1, v_2)\> \right].
\end{aligned}
\end{equation}
We then choose $f(x)$ to be the the fundamental solution $p^h(t, x, y_0)$, so
$P_t^hf(x_0, y_0)=p(2t, x_0, y_0)$. In particular,
\begin{equation}\begin{aligned}
\f{  \Hess p^{h} (2t, x_0, y_0)(v_1, v_2)} { p^{h} (2t, x_0, y_0)}
=& \EE\left[\f {  p(t, x_t, y_0)} { p(2t, x_0, y_0)} N_{ t}\right]
+ \EE \left[ \f {  p(t, x_t, y_0)} { p(2t, x_0, y_0)} \f 2 t \int_0^{t/2} \<X(x_s) dB_s,  W_{s}^{(2)}(v_1, v_2)\> \right].
\end{aligned}
\end{equation}
The right hand side can then be estimated.  Since $|W_t|$ is bounded by a deterministic function
 (when $\rho^h$ is bounded above), the first term of the right hand side  is easier to estimate. Let us work with the second term, 
\begin{equation*}\begin{aligned}
&\f {1} {p(2t, x_0, y_0)} \EE \left( p(t, x_t, y_0)\f 2 t\int_0^{t/2} \<X(x_s) dB_s,  W_{s}^{(2)}(v_1, v_2)\> \right)\\
 &\le \f 2 t \EE  \left(\f{p(t, x_t, y_0) }{p(2t, x_0, y_0)}   \log \f{p(t, x_t, y_0) }{p(2t, x_0, y_0)}   \right) 
 +\f 2 t \log \EE \left ( \exp\left(\int_0^{t/2} \<X_i(x_s) dB_s,  W_{s}^{(2)}(v_1, v_2)\>\right)\right)\\
 \le & \f 2 t \sup_{y\in M}   \log \f{p(t, y, y_0) }{p(2t, x_0, y_0)}  
 +\f 2 t \log \EE \left ( \exp\left(\int_0^{t/2} \<X_i(x_s) dB_s,  W_{s}^{(2)}(v_1, v_2)\>\right)\right).
\end{aligned}
\end{equation*}
This can then be refined by  heat kernel estimates and by estimates on $\EE\left( e^{ \left|W_{s}^{(2)}(v_1, v_2)\right|^2}\right)$, we illustrate the latter below. The other terms can be treated similarly.

\begin{lemma}\label{exponential-estimates}
Suppose that $|\rho^h| \le K$,  $\|\Rc_x\|\le \|\Rc\|_{\infty}$,  and $\|\Theta^h\|^2\le c+ \delta r^2$ for $\delta$ sufficiently small.
 Set $C_1(T,0)=1$,
$$C_1(T,K)= \sup_{0< s\le 3KT}\f 1{s}{(e^s-1)}, \quad \alpha_2(T, K, \|\Rc\|_\infty)=\f  1{49 n^2\|\Rc\|_\infty^2 C_1(T,K)}.$$
Then 
there exists a universal constant $c$ such that for unit vectors $v_1, v_2\in T_{x_0}M$,
and for any $\alpha \le  \alpha_2(T, K, \|\Rc\|_\infty)$, 
\begin{eqnarray*} \EE \exp\left(\alpha  \gamma |W_t^{(2)}(v_1, v_2)|^2\right)
& \le& c e^{2\alpha \gamma} 
\sqrt{\EE \exp\left(4t \;\gamma \alpha   \int_0^t e^{3 Ks} \|\Theta^h\|^2_{x_s} ds\right)}
\\
&\le &c e^{ \f  {2 \gamma}{49 n^2\|\Rc\|_\infty^2}} \sqrt{ \EE \exp\left(\f{4t\gamma}{49n^2\|\Rc\|_\infty^2C_1(t,K)}  \int_0^t e^{3 Ks} \|\Theta^h\|^2_{x_s} ds\right)}.\end{eqnarray*}
\end{lemma}

At this stage we must choose an optimal  condition on the growth of $|\Theta^h|$ at the infinity  and estimate the exponential integrability of the radial distance function. The condition we imposed is linear growth, with
the linear part sufficiently small (or we may compensate the size of the linear part by taking $t$ in a small interval $[0, t_0]$).
With these estimates we conclude this paper, and invite  the interested reader to  consult \cite{Hessian}
for technicalities and further results. There we also studied the class of  manifolds with a pole 
 using semi-classical bridges. The use of semi-classical bridge for derivatives estimates is novel.
 See also \cite{Li-ibp-sc}.
 Finally we pose the following open question. We know that a bound of the form $e^{-Kt}$ on $P_t$ characterises the lower boundedness of a Ricci curvature. Elworthy asked me whether I can use equation (\ref{basic-formula}) and obtain some characterisation for manifolds. Let me make  precise a question here.

{\bf Open Problem.} 
For $f\in BC^2$, by  Lemma \ref{second-diff-lemma-1},
$$\left| \Hess (P^h_tf)\right| _{L_\infty}
\le | \nabla df|_\infty \;  \EE  \left( e^{2\int_0^t\rho^h(x_s) ds}\right)
+|df|_{L_\infty}\; \EE  \left| W^{(2)}_t\right|.
$$
If the Ricci curvature is bounded from below by $K$, we have
$$\left| \Hess (P^h_tf)\right| _{L_\infty}
\le | \nabla df|_\infty e^{-2Kt}
+|df|_{L_\infty}\; \EE  \left| W^{(2)}_t\right|.
$$
Can we characterise the class of complete Riemannian manifolds, among those whose Ricci curvature is bounded from below by $K$ and  whose sectional curvature and  symmetrised tensor $\Theta^h$ are bounded?


\def\cprime{$'$} \def\cprime{$'$} \def\cprime{$'$} \def\cprime{$'$}

\end{document}